\documentclass[a4paper,11pt]{article}
\usepackage{amsfonts,amssymb}

\newtheorem{lemma}{Lemma}

\newtheorem{prop}{Proposition}

\newtheorem{cor}{Corollary}

\newtheorem{thm}{Theorem}

\def\pf{{\it Proof. }}

\def\CC{{\mathbb C}}
\def\ZZ{{\mathbb Z}}

\def\calA{{\mathcal A}}
\def\calC{{\mathcal C}}
\def\calI{{\mathcal I}}
\def\calS{{\mathcal S}}
\def\calT{{\mathcal T}}

\def\longra{\longrightarrow}
\def\da{\downarrow}
\def\implies{\Rightarrow}

\def\vep{\varepsilon}

\begin{document} 

\begin{center}{\Large\bf 
Arrangements associated to chordal graphs\\\,  
and limits of colored braid groups}
\end{center}

\begin{center} 
{\large\bf Fr\'ed\'eric Chapoton and Patrick Polo}
\end{center}

\bigskip 
\begin{abstract} 
Let $G$ be a chordal graph, $X(G)$ the complement of
the associated complex arrangement and $\Gamma(G)$ the fundamental
group of $X(G)$. 
We show that $\Gamma(G)$ is a limit of colored braid groups over the poset 
of simplices of $G$. 
When $G = G_T$ is the comparability graph associated with a rooted tree $T$, 
a case recently investigated by the first author, the result takes the following
very simple form: 
$\Gamma(G_T)$ is a limit over $T$ of colored braid groups.
\end{abstract}


\bigskip 
\section*{Introduction}

In \cite{Ch}, the first author investigated certain properties of the 
arrangements of hyperplanes associated to rooted trees. 
As pointed out by Stanley, these 
arrangements are special cases of arrangements associated to chordal 
graphs. Let $G$ be such a graph, let $X(G)$ denote the complement of
the associated complex arrangement and let $\Gamma(G)$ be the fundamental
group of $X(G)$. 

In this paper, we show that $\Gamma(G)$ can be described, in a compact manner, 
as a limit of colored braid groups over the poset $\calS(G)$ 
of simplices of $G$ 
(that is, the subsets of the set of vertices for which the induced subgraph
is complete).  
In the special case where $G = G_T$ is the comparability graph associated
with a rooted tree $T$, the maximal simplices of $G_T$ are, of course, parametrized
 by the leaves of $T$ and the result can be reformulated as follows: 
 $\Gamma(T) := \Gamma(G_T)$ is a limit over $T$ of colored braid groups.
 
\medskip The plan of the article is as follows. In Section 1, we
recall several facts about chordal graphs, graphical arrangements, fibrations,
pull-backs and limits. After these preliminaries, we state and prove, in
Section 2, the main theorem. Finally, in Section 3, we consider the
case of rooted trees, where the results take a simpler and slightly more precise
form.

\section{Recollections}

\subsection{Chordal graphs}\label{chord}

Let $G = (V,E)$ be a finite graph, where $V$ is the set of vertices and
$E$ the set of edges. A cycle $C$ in $G$, of length $k\geq 3$,  is a
sequence $v_1,\dots,v_k$ of pairwise distinct vertices such that 
$\{v_k,v_1\}$ and $\{v_i,v_{i+1}\}$, for $i=1,\dots,k-1$, are edges. 
A {\it chord\,} of $C$ is an edge $\{v_i,v_j\}$ distinct from the previous
 ones, {\it i.e.}, such that $v_i$ and $v_j$ are not 
consecutive vertices in the cycle. 

One says that $G$ is {\it chordal\,} (or a rigid circuit graph) if every cycle of
length $\geq 4$ admits at least one chord. 

\smallskip 
Any subset $S$ of $V$ determines an induced subgraph, whose edges are the
edges of $G$ connecting two elements of $S$. 
By abuse of language, 
we will most of the time make no distinction between a subset of $V$
and the corresponding subgraph. 
{F}rom the definition above it is clear that if $G$ is chordal then so is
the subgraph induced by any subset $V'$ of $V$. 

\medskip On the other hand, one says that a subset $S$ of $V$ 
is a simplex
 of $G$ if the induced subgraph is complete, that is, if $\{s,s'\}$ is
 an edge for all $s\not= s'$ in $S$. 
 
\subsection{Simplicial vertices and PEO's}\label{PEO} 

{F}or any vertex $v$, let $N(v)$ denote the set of its neighbours, that is,
the vertices connected to $v$. One says that $v$ is a 
{\it simplicial vertex\,} if $N(v)$ is a simplex. In this case, 
$N(v)\cup\{v\}$ is a simplex, too. 

{F}ollowing \cite{BP}, we say that an ordering $v_1,v_2,\dots,v_n$, 
where $n = \vert V\vert$, 
of the vertices of $G$ is a perfect elimination ordering (PEO in short)
if $v_i$ is a simplicial vertex of $G\setminus \{v_1,\dots,v_{i-1}\}$, for all 
$i = 1,\dots,n$. (This is called a vertex elimination order in \cite{ER};
we prefer the above terminology, which we find more suggestive.) 

Of course, in a complete graph any vertex is simplicial and any ordering
of the vertices is a PEO.  

{F}or future reference, let us record the following lemma and corollary.

\begin{lemma} {\rm 1)} A vertex is simplicial if and only if it belongs to a unique
maximal simplex of $G$. 

{\rm 2)} Suppose that $G$ is chordal but not complete. Then $G$ admits
at least two nonadjacent simplicial vertices. In particular, if $S$ is
a simplex of $G$, then $G\setminus S$ contains a simplicial vertex of $G$.
\end{lemma} 

The first part is an almost immediate consequence of the definitions, 
see e.g. \cite[Lemma 3.1]{BP}. For the second part, we refer to 
\cite[Lemma 2.2]{BP} or \cite[Lemma 4.2]{Go}.

\medskip Since any induced subgraph of a chordal graph is still chordal,
this implies immediately the 

\begin{cor} Suppose that $G$ is chordal and let $S$ be any simplex in 
$G$. Then there exists a PEO $v_1,\dots,v_n$ of $G$ such that the vertices 
not in $S$ appear first, {\it i.e.}, such that 
$\{v_1,\dots,v_k\} = V\setminus S$, where $k = \vert V\setminus S\vert$.
\end{cor} 

It turns out that, conversely, the existence of a PEO implies
that $G$ is chordal, see \cite[Th. 4.1]{Go} or \cite[Th. 2.3]{BP}. We will not
need this fact. 

\subsection{Fibrations}\label{fibrations}

We gather in this subsection several concepts and results 
from homotopy theory that will be used in the sequel. 

By a space we mean a pointed, Hausdorff topological space and 
maps are base-point preserving continuous maps. 
Let $B$ be a pathwise-connected space, with base point $b_0$. 
A continuous map $f : E\to B$ is called a Hurewicz fibration, 
or simply fibration, if it has the 
homotopy lifting property with respect to any space. 

{F}irstly, the composite of two fibrations is again a fibration. 
Secondly, if $f$ is a locally trivial fibration, it is a Hurewicz 
fibration, at least when the base space $B$ is paracompact, e.g. 
metrisable. Thirdly, if $f$ is a fibration, all fibers 
$f^{-1}(b)$ are homotopy equivalent to $F := f^{-1}(b_0)$. 
Fourthly, assuming for simplicity that $F$ is pathwise connected, 
there is a long exact 
sequence of homotopy groups 
$$
\cdots\to \pi_2(B)\to \pi_1(F)\to \pi_1(E)\to \pi_1(B)\to 1. \eqno (\dag)
$$
Moreover, if $f$ is split, {\it i.e.}, admits a section, then $(\dag)$ 
breaks into split short exact sequences. 
For all this see, for example, \cite[Chap.7]{Se} 
or \cite[\S\S I.7 \& IV.8]{Wh}.

\medskip We will also need the following functorial property of the 
exact homotopy sequence of a fibration, see \cite[Chap.6]{Hi}. 
Consider a (strict) pull-back diagram 
$$
\begin{array}{rcl}
W  & \stackrel{\tilde{g}}{\longra} & X\\
{\scriptstyle \tilde{f}} \da& & \da {\scriptstyle f}\\
Y & \stackrel{g}{\longra}& B.
\end{array}
$$
By this, we mean that $W$ identifies under $(\tilde{f},\tilde{g})$  
with the closed subset of $Y\times X$ consisting of 
those pairs $(y,x)$ such that $g(y) = f(x)$. 
Suppose further that $f$ is a fibration, and let $F := f^{-1}(b_0)$ 
and 
$\widetilde{F} := \tilde{f}{}^{-1}(y_0)$. 
 
Then $\tilde{f}$ is also a fibration and, under the identification above, 
the map $\widetilde{F} \to F$ induced by $\tilde{g}$ is the identity. 
For simplicity, assume further that $F$ is pathwise connected. Then, 
one has a commutative ladder 
$$
\begin{array}{rcccl}
\cdots\to\pi_2(Y) & \to \pi_1(F)\to & \pi_1(W) & \stackrel{\tilde{f}}{\longra}&
\pi_1(Y) \to 1\\
{\scriptstyle g}\da\phantom{Y}& {\scriptstyle 1}\da& {\scriptstyle \tilde{g}}\da& &
\phantom{\pi}{\scriptstyle g}\da\\
\cdots\to\pi_2(B) & \to \pi_1(F)\to & \pi_1(X) & \stackrel{f}{\longra}&
\pi_1(B) \to 1.
\end{array}
$$

{F}rom this, one easily obtains the following 

\begin{cor} If $\pi_1(F)\to \pi_1(X)$ is injective (e.g., 
if $\pi_2(B) = \{1\}$, or if $f$ admits a section), then so is  
$\pi_1(F)\to \pi_1(W)$. In this case, $\pi_1(W)$ is the pull-back 
of the diagram 
$$
\pi_1(X)\stackrel{f}{\longra}\pi_1(B)\stackrel{g}{\longleftarrow}
\pi_1(Y).
$$
\end{cor}

\subsection{Graphical arrangements and 
locally trivial fibrations}\label{simplicial-lemma}

Let $G = (V,E)$ be a finite graph. 
It determines an arrangement $\calA(G)$ in $\CC^G$, whose hyperplanes are
given by the linear forms $x_i - x_j$, when $\{i,j\}\in E$.  
Let us denote by $X(G)$ the complement of $\calA(G)$ and by 
$\Gamma(G)$ the fundamental group of $X(G)$.  

\medskip Let $v\in V$. Let us denote by $p_v$ the restriction 
to $X(G)$ of the natural projection 
$\CC^G\to \CC^{G\setminus v}$. Clearly, it maps 
$X(G)$ onto $X(G\setminus v)$ and, for every 
$\mathbf{z} = (z_u)_{u\in G\setminus v}\in X(G\setminus v)$, 
one has  
$$
p_v^{-1}(\mathbf{z}) = \{x\in \CC\mid x\not= z_u, \quad\forall u\in N(v) \}.
$$
Set $\nu(\mathbf{z}) = \vert \{z_u \mid u\in N(v)\}\vert$. 
Of course, if $p_v$ is a locally trivial fibration then 
$\nu(\mathbf{z}) = \vert N(v)\vert$ for all 
$\mathbf{z} \in X(G\setminus v)$. It is well-known that the converse 
holds, see for example \cite[Lemma 1]{Ar} or \cite[Prop.2]{Br}. 
This fact is also used implicitly in the proof of 
\cite[Th. 2.9]{Te86}, where a considerably more general 
result is proved (see also \cite[Th. 5.111]{OT}). 
Since the details are not easily found in the literature 
(see, however, the proof of \cite[Th.1]{FN}), 
we recall the argument for the reader's (and for our own)
convenience. 

\begin{lemma} 
Let $v$ be a simplicial vertex with $N$ neighbours. Then the projection 
$p_v : X(G)\to X(G\setminus v)$ is a split, locally trivial smooth fibration, with 
fiber $\CC $ minus $N$ points. 
\end{lemma} 

\pf Let us choose a numbering 
$u_1,\dots, u_n$ of the vertices of $G\setminus v$ 
such that $N(v) = \{u_1,\dots,u_N\}$, and denote 
the corresponding coordinates $z_{u_i}$ on 
$X(G\setminus v)$ simply by $z_i$. 

Let $\mathbf{z}^0\in X(G\setminus v)$. 
Since $v$ is a simplicial vertex, one has $z_i^0\not= z_j^0$ for 
$i\not= j$ in $[1,N]$.
Pick $\vep >0$ such that the open disks $D(z^0_i,2\vep)$ are disjoint,
and let 
$U$ denote the open subset of $X(G\setminus v)$ defined by $\vert z_i -z_i^0\vert
< \vep$, for $i= 1,\dots,N$. 

Let $\theta : \CC\to [0,1]$ be a smooth function such that 
$\theta(z) = 1$ if $\vert z\vert \leq \vep$ and 
$\theta(z) = 0$ if $\vert z\vert \geq 3\vep/2$. 
The map $V : U\times \CC\to \CC$ given by 
$$
V(\mathbf{z},x) = \sum_{i=1}^N \theta(x-z_i^0) (z_i-z_i^0)
$$
is clearly $C^\infty$. For each $\mathbf{z}$, it defines a vector field 
$V_{\mathbf{z}}$ on $\CC$, which depends smoothly on the parameter $\mathbf{z}$,
 vanishes outside the union of the disks $D(z_i^0,3\vep/2)$, 
 and coincides on each $\overline{D}(z_i^0, \vep)$ with the constant vector field 
 $z_i - z_i^0$. 
  Since $V_{\mathbf{z}}$ has compact support, the 
 flow $\psi(\mathbf{z},t,x)$ corresponding to the initial condition 
 $\psi(\mathbf{z},0,x) = x$, is defined for all values of $t$. 
 
 Set $\phi(\mathbf{z},x) = \psi(\mathbf{z},1,x)$ and 
 $\phi_{\mathbf{z}}(x) = \phi(\mathbf{z},x)$. Then $\phi$ is 
 $C^\infty$ and it is not difficult to see that, for every 
 $\mathbf{z}$, $\phi_{\mathbf{z}}$ is a diffeomorphism of $\CC$ 
 which equals the identity outside the union of the disks $D(z_i^0,3\vep/2)$
 and maps each $z_i^0$ to $z_i$. 
Moreover, $\phi_{\mathbf{z}}^{-1}(x)$ is nothing but 
$\psi(\mathbf{z},-1,x)$ and hence the map 
$\mu : (z,x)\mapsto \phi_{\mathbf{z}}^{-1}(x)$ is $C^\infty$, too. 
Therefore, the maps  
$$
f : U\times \left(\CC\setminus\{z^0_1,\dots,z^0_N\}\right) \rightleftarrows\, 
p_v^{-1}(U) : g
$$
given by 
$f(\mathbf{z},x) = (\mathbf{z},\phi(\mathbf{z},x))$
 and $g(\mathbf{z},x) = (\mathbf{z},\mu(\mathbf{z},x))$
 are inverse diffeomorphisms. This proves that $p_v$ is a smooth locally
 trivial fibration.  
 
 \medskip Moreover, by \cite[Lemma 1]{Ar} or \cite[Prop. 2.4]{FR85}, 
 $p_v$ admits a section $\sigma_v$. Indeed, following \cite{FR85}, 
 one may take 
 $\sigma_v(\mathbf{z}) = \left(\mathbf{z}, 
 (1+\sum_{i=1}^N \vert z_i\vert{}^2){}^{1/2}\right)$.
 This completes the proof 
 of the lemma. 
 
 \bigskip As a well-known consequence of the lemma (see \cite{FR85}),
one obtains that the long exact homotopy sequence of the locally trivial fibration 
 $p_v$ gives isomorphisms 
 $\pi_i(X(G))\cong \pi_i(X(G\setminus v))$
 for $i\geq 2$, along with a split exact sequence 
 $$
 1\to F_N \to \Gamma(G)\to \Gamma(G\setminus v)\to 1,
 $$
 where $F_N$ is the free group on $N$ generators, $N = \vert N(v)\vert$.

\subsection{Chordal graphs and fibrations}\label{fib-simplexe}

Return to $X(G)$, the complement of the arrangement 
associated to our graph $G$. By the previous paragraph, 
if $v$ is a simplicial vertex of $G$, then 
$X(G)\to X(G\setminus v)$ is a locally trivial fibration 
with a section.   

Assume now that $G$ is chordal. 
Of course then, by iterated use of the previous lemma
one obtains that $\calA(G)$ is of fiber-type in the sense of Falk and Randell
\cite{FR85} and hence that $X(G)$ is a $K(\pi,1)$ space. 

Moreover, since a composite of fibrations 
(resp. sections thereof)
is a fibration (resp. a section thereof), 
one obtains, using Corollary \ref{PEO}, the following corollary, 
which will prove useful later.

\begin{cor} Let $S$ be a subset of $V$ such that there exists a PEO 
of $G$ beginning with the elements of $V\setminus S$. 
Then $X(G)\to X(S)$ 
is a split fibration, whose 
fiber and base are $K(\pi,1)$ spaces. 
This is the case, in particular, if $S$ is a simplex of $\,G$.
\end{cor} 

\medskip\noindent{\bf Remark.} It follows from \cite[Th.\,2.2 \& Cor.\,2.4]{Pa} 
that the projection $X(G) \to X(S)$ is actually a locally trivial
smooth fibration. 
This is a general result about the projection corresponding
to any modular element in the lattice of an arbitrary hyperplane arrangement, 
whose proof relies on Thom's first isotopy lemma. We will not need
this more general result. 

\subsection{Colored braid groups}\label{braid}

Let $\calI$ denote the category whose objects are finite sets and whose
morphisms are the injective maps. For any finite set $I$, let $C_I$ denote
the complete graph on $I$ and let 
$$
P(I) := \Gamma(C_I) = \pi_1(X(C_I)).
$$
Thus, $P(I)$ is isomorphic to the braid group on $n$ colored strands,
where $n = \vert I\vert$, see \cite{FoN} or \cite[\S1.4]{Bi}. 

\begin{prop} The assignment $I\mapsto P(I)$ is a contravariant functor
from $\calI$ to the category of groups. 
\end{prop} 

\pf For every subset $J$ of $I$, one has the natural projection 
$\phi_{JI} : X(C_I)\to X(C_J)$ and the associated group homomorphism 
$\rho_{JI} := \pi_1(\phi_{JI})$, from $P(I)$ to $P(J)$. 
Since, clearly, $\phi_{KJ}\circ\phi_{JI} = \phi_{KI}$, 
when $K\subset J\subset I$, the result follows from the functoriality 
of $\pi_1$.

\subsection{Limits}\label{limits}

Let $\calS$ be a finite poset and let $F$ be a contravariant functor 
from $\calS$ to the category of groups. That is, one is given 
a morphism $\rho_{ts} : P(s)\to P(t)$ whenever $t\leq s$, 
such that $\rho_{ut}\rho_{ts} = \rho_{us}$ if $u\leq t\leq s$. 

Then, the limit of $F$ over $\calS$ is the following subgroup 
of $\prod_{s\in\calS} F(s)$ : 
$$
\lim_{\calS} F := \left\{ (g_s)_{s\in\calS}\in \prod_{s\in\calS} F(s)\mid
\rho_{ts}(g_s) = g_t, \quad\forall t\leq s \right\}.
$$

{\bf Remarks.} 1) This is, of course, a special case of 
the general notion of limit,  
see for instance \cite[Chap.V]{ML} or \cite[Appendix A.5]{We}. 
The above definition is sufficient 
for our purposes. 

2) Denoting by $M(\calS)$ the set of maximal elements of $\calS$,
$\lim_{\calS} F$ may also be regarded as a subgroup of 
the (smaller) product $\prod_{s\in M(\calS)} F(s)$.

\medskip 
{F}or any $s\in\calS$, set 
$D(s) := \{t  \mid t\leq s \}$. 
Further, if $s$ is a maximal element of $\calS$, define its cone   
$$
\calC(s) = \{t \leq s \mid \forall u\in \calS,\, t\leq u\Rightarrow
u\leq s\}. 
$$
Then, the following lemma is straightforward and its proof 
left to the reader. 

\begin{lemma} Let $s$ be a maximal element of $\calS$ and let 
$\calT$ be a subset of $\calC(s)$ containing $s$. 
Then 
$\lim_{\calS} F$ is isomorphic to the pull-back of 
the diagram 
$$ 
F(s)\to \lim_{D(s)\setminus \calT} F \leftarrow \lim_{\calS\setminus \calT} F.
$$
\end{lemma}

\section{Limits of colored braid groups}

\subsection{The poset of simplices and the theorem}\label{main}  

Again, we return to our chordal graph $G = (V,E)$ and 
recall that $\Gamma(G) := \pi_1(X(G))$.   
Let us denote by $\calS(G)$ the poset of non-empty simplices of $G$,
ordered by inclusion. Observe that $\calS(G)$ is a meet-semilattice. 

We will regard $\calS(G)$ as a subposet of the lattice 
of subsets of $V$. That is, we identify each 
simplex of $G$ with the underlying subset of $V$. Then, 
recall from \ref{braid} the functor 
$P : \calI \to$ $\{$Groups$\}$. 
We can now state our main result. 

\begin{thm} Let $G$ be a chordal graph. Then $\Gamma(G)$ is the limit 
over $\calS(G)$ of the functor $P$, that is, 
$$
\Gamma(G) = \lim_{\calS(G)} P = \left\{ (\gamma_S)_{S\in\calS(G)}\in \prod_{S\in\calS(G)} P(S)\mid
\rho_{TS}(\gamma_S) = \gamma_T, \quad\forall T\leq S \right\}.
$$
\end{thm}

The theorem is of course true if $G$ is complete or, more generally, a
disjoint union of complete graphs. The proof will be 
by induction on the number of vertices of $G$.

\subsection{A pull-back diagram of spaces}

Let $s_0$ be a simplicial vertex of $G$ and let $S = N(s_0)\cup \{s_0\}$;
it is the unique maximal
simplex containing $s_0$ (see Lemma \ref{PEO}). Let us introduce 
$$
S^0 := \{s\in S\mid \forall t\in V, \quad \{s,t\}\in E\implies t\in S\}.
$$
Note that $S^0$ is not empty since $s_0\in S^0$. 

With obvious notation, there is a pull-back diagram 
$$
\begin{array}{ccc}
\CC^G & \longra & \CC^{G\setminus S^0}\\
\da   &         & \da \\
\CC^S & \longra & \CC^{S\setminus S^0}.
\end{array}
$$

\begin{lemma} The previous diagram induces a pull-back diagram 
$$
\begin{array}{ccc}
X(G) & \longra & X(G\setminus S^0)\\
\da   &         & \da {\scriptstyle f} \\
X(S) & \stackrel{\phi}{\longra} & X(S\setminus S^0).
\end{array}\eqno (*)
$$
\end{lemma}

\pf  One observes first that the lemma is implied by the following claim:
if $\{i,j\}$ is an edge
of $G$ then $i,j\in S$ or $i,j\in G\setminus S^0$. 

To prove this last assertion, suppose that $\{i,j\}$ is an edge
of $G$. It is of the desired type if $i,j\in S$ or $i,j\not\in S^0$. 
So, the only case to worry about is the case where $i\in S^0$ and $j\not\in S$. 
But then $\{i,j\}$ is not an edge, by the definition
of $S^0$. This proves the assertion and, hence, the lemma.

\subsection{Completion of the proof}
Thus, $(*)$ is a pullback diagram. Moreover, since $S\setminus S^0$ is a simplex 
of $G\setminus S^0$, it follows from Corollary \ref{fib-simplexe} that
 $f$ is a split fibration. 
(We will see below that, for a slightly different, though similar, reason, 
$\phi$ is also a split fibration, but we do not need this for the moment).

Therefore, one deduces 
from Corollary \ref{fibrations} that $\Gamma(G)$ is the pull-back 
of the diagram 
$$
\Gamma(S) \longra  \Gamma(S\setminus S^0)\longleftarrow \Gamma(G\setminus S^0).
$$

{F}urthermore, $\Gamma(S) = P(S)$ and  
$\Gamma(G\setminus S^0)\cong \lim_{\calS(G\setminus S^0)} P$, 
by the inductive hypothesis. Since $S\setminus S^0$ is a simplex, 
one has, moreover, 
$\Gamma(S\setminus S^0) = P(S\setminus S^0) \cong \lim_{\calS(S\setminus S^0)} P$.
Therefore, one obtains a pull-back diagram of groups:
$$
\begin{array}{ccc}
\Gamma(G)&\longra & \lim\limits_{\calS(G\setminus S^0)} P\\
\da& & \da \\
P(S)& \longra & \lim\limits_{\calS(S\setminus S^0)} P.
\end{array} \eqno (**)
$$

\begin{lemma} {\rm 1)} Every simplex of $G$ which meets $S^0$ is contained
in $S$. 

\noindent {\rm 2)} Every element of $S^0$ is a simplicial vertex of $G$.
\end{lemma} 

\pf Let $s\in S^0$. By the very definition of $S^0$, no neighbour of $s$
belongs to $G\setminus S$. This implies the first assertion, and the
second follows immediately. 

\medskip 
Denote by $\calT$ the set of simplices of $G$ which meet $S^0$ and 
recall from \ref{limits} the definition of $\calC(S)$, the cone of 
$S$. Then, the lemma shows that $\calT$ is contained in $\calC(S)$. 
Moreover, one has, clearly, 
$\calS(G\setminus S^0) = \calS(G)\setminus \calT$ and 
$\calS(S\setminus S^0) = \calS(S)\setminus \calT$. 
Therefore, using Lemma \ref{limits}, one deduces from $(**)$ 
that 
$$
\Gamma(G) \cong \lim_{\calS(G)} P.
$$
This completes the proof of the theorem. 

\medskip {\bf Remark.} One deduces from the second assertion of the lemma,
coupled with Corollary \ref{fib-simplexe}, that in fact all maps 
in the diagram $(*)$ are split fibrations.

\section{The case of rooted trees}

In the case of rooted trees, the results take an even more precise 
form, which we describe in this final section. 

\subsection{The theorem for rooted trees} 
Let $T$ be a rooted tree, let $S$ be its set of vertices, and let 
$r$ be the distinguished vertex called the root.  
For each vertex $s$, let $T_s$ be 
the linear tree joining $s$ to $r$. 
Let $\leq$ denote the partial order on $S$ induced by $T$, 
that is, $s \leq t$ iff $T_s\subseteq T_t$; 
 the maximal elements are exactly the leaves of $T$. 

\smallskip Recall the notation of \ref{braid} and, for each $s$, 
set $P(s) := P(T_s)$. 
If $s'\leq s$, then $T_{s'}\subseteq T_s$ and hence one has 
a morphism $\rho_{s's} : P(s) \to P(s')$. Therefore, $P$ 
may be regarded as a functor on the poset $T$. 

\medskip
 Let us denote by $\calA_T$ the arrangement associated to $T$ 
(see \cite{Ch}) and by $\Gamma_T$ the fundamental group of the complement. 
Thus, if $G_T$ denotes the comparability graph of $T$, one 
has $\calA_T = \calA(G_T)$ and $\Gamma_T = \Gamma(G_T)$.

\smallskip 
Clearly, the leaves of $T$ are simplicial vertices of $G_T$ 
and any ordering $v_1,v_2,\dots$ of the vertices such that 
each $v_i$ is a leaf of $T\setminus \{v_1,\dots,v_{i-1}\}$ 
is a PEO. In particular, $G_T$ is chordal ({\it cf.} \cite[Lemma 2]{Ch}).

Further, 
the maximal simplices of $G_T$ are the linear trees $T_f$, 
for $f$ a leaf of $T$. 
Moreover, any (finite) intersection of such 
trees is equal to a tree $T_s$, for some $s\in S$. 
From this, one easily
deduces that Theorem \ref{main} can be reformulated as follows. 

\begin{thm} One has $\Gamma_T \cong \lim_T P$. 
\end{thm}

\subsection{The structure of iterated semi-direct product}

Moreover, in the case of rooted trees, the structure of $\Gamma_T$ as 
an iterated semi-direct product of free groups 
(which is a general feature of fundamental groups of fiber-type 
arrangements, see \cite{FR85}), can be described rather precisely, 
as follows.

Recall that each fibration $X(T_s)\to X(T_{s'})$ admits 
a section, say $\tau_{ss'}$. Therefore, 
each  $P(s)$ may be identified to a semi-direct product 
$$
P(s) \cong K_{s's}\ltimes \tau_{ss'}(P(s')), 
$$
where $K_{s's} = \ker \rho_{s's}$. Then, one obtains easily the 
following proposition, whose proof is left to the reader. 
Define the {\em height} of any vertex $s$ by 
$h(s) := \vert T_s\vert -1$ and let  
$n=n(T)$ be the supremum of the $h(f)$, for $f$ a leaf of $T$. 

\begin{prop} There is an isomorphism 
$$
\Gamma_T\cong \prod_{h(s_1)=1} \ZZ(s_1)\ltimes\left(\prod_{h(s_2)=2} F_2(s_2)\ltimes
\cdots 
\ltimes\left(\prod_{h(s_n)=n} F_n(s_n)\right)
\right)
$$
(iterated semi-direct product, for $i=1,\dots,n$, of 
direct products, indexed by the vertices of height $i$,  
 of copies of the free group $F_i$), 
 such that, for every $t\in S$, 
the projection $\Gamma_T\to P(t)$ identifies  
 with the natural projection from $\Gamma_T$ to 
$$
\ZZ(s_1(t))\ltimes \Bigl( F_2(s_2(t))\ltimes \cdots \ltimes F_h(s_h(t))\Bigr),
$$
where  $h = h(t)$ and $(r, s_1(t),\dots,s_h(t))$ is the ordered chain of 
vertices $\leq t$. 
\end{prop}

\bigskip{\small
\begin{flushright}
{F}r\'ed\'eric Chapoton\\
LACIM\\
Universit\'e du Qu\'ebec \`a Montr\'eal\\
Montr\'eal Qu\'ebec H3C 3P8\\
Canada\\
chapoton@math.uqam.ca
\end{flushright} 

\begin{flushright}
Patrick Polo\\
CNRS \& Universit\'e Paris 13\\ Institut Galil\'ee (LAGA)\\ 
{F}-93430 Villetaneuse\\
{F}rance\\
polo@math.univ-paris13.fr
\end{flushright}
}

\end{document}